\documentclass{article}
\usepackage{amsmath}
\usepackage{amssymb}
\usepackage{amsfonts}

\newtheorem{theorem}{Theorem}

\newtheorem{corollary}[theorem]{Corollaire}

\newtheorem{lemme}[theorem]{Lemme}

\newtheorem{remarque}[theorem]{Remarque}

\newenvironment{proof}[1][Proof]{\noindent\textbf{#1.} }{\ \rule{0.5em}{0.5em}}
\newenvironment{preuve}[1][Preuve]{\noindent\textbf{#1.} }{\ \rule{0.5em}{0.5em}}
\newenvironment{Preuve du lemme}[1][Preuve du lemme]{\noindent\textbf{#1.} }{\ \rule{0.5em}{0.5em}}
\include {mak}
\input{tcilatex}
\begin{document}

\title{ G\'{e}n\'{e}ralisation des congruences de Wolstenholme et de Morley}
\author{ Farid Bencherif et Rachid Boumahdi \\
\ \ \ \ \ \ \ \ \ \ \ \ \ \ \ \ \ \ \ \ \ \ \ \ \ \ \ \ \ \ \ \ \ \ \ \ \ \
\ \ \ \ \ \ \ \ \ \ \ \\
{\small Laboratoire LA3C USTHB, Fac. Math. P.B. 32, El Alia, 16111, Algiers,
Algeria.}\\
{\small \ fbencherif@gmail.com, r\_boumehdi@esi.dz }\\
}
\date{}
\maketitle

\begin{abstract}
Dans cet article, nous prouvons que pour tout nombre premier impair $p$ et
pour tout $p$-entier $\alpha $, on a 
\begin{equation*}
\binom{\alpha p-1}{p-1}\equiv 1-\alpha (\alpha -1)(\alpha ^{2}-\alpha
-1)p\sum_{k=1}^{p-1}\frac{1}{k}+\alpha ^{2}(\alpha -1)^{2}p^{2}\sum_{1\leq
i<j\leq p-1}\frac{1}{ij}\text{ }\func{mod}(p^{m}).
\end{equation*}%
o\`{u} $m=7$ si $p\neq 7$ et $m=6$ si $p=7$.

Cette congruence g\'{e}n\'{e}ralise les congruences de Wolstenholme, Morley,
Glaisher, Carlitz, McIntosh, Tauraso et Me\v{s}trovi\'{c}. Elle permet aussi
de retrouver simplement des congruences d\^{u}es \`{a} Glaisher, Carlitz et
Zhao. 
\begin{equation*}
\end{equation*}


\end{abstract}

\section{Introduction}

Dans tout ce qui suit $p$ d\'{e}signe un nombre premier impair.

En 1819, Babbage \cite{bab} prouve que pour tout nombre premier $p\geq 3$,
on a%
\begin{equation*}
\binom{2p-1}{p-1}\equiv 1\text{ }\func{mod}(p^{2}).
\end{equation*}

En 1862, Wolstenholme \cite{wol} et $[$\cite{har}, p. 89$]$ prouve que pour
tout nombre premier $p\geq 5$, on a les deux congruences suivantes:%
\begin{equation}
\binom{2p-1}{p-1}\equiv 1\text{ }\func{mod}(p^{3}),  \label{rel70}
\end{equation}%
\begin{equation*}
\sum_{k=1}^{p-1}\frac{1}{k}\equiv 0\text{ }\func{mod}(p^{2}).
\end{equation*}%
En 1895, Morley \cite{mor} prouve que 
\begin{equation}
(-1)^{\frac{p-1}{2}}\binom{p-1}{\frac{p-1}{2}}\equiv 4^{p-1}\text{ }\func{mod%
}(p^{3})\text{ \ \ \ \ \ \ }p\geq 5.  \label{rel71}
\end{equation}%
En1900, Glaisher $\left[ \text{\cite{gla1}},\text{ p. 21}\right] $, $\left[ 
\text{\cite{gla2}},\text{ p. 323}\right] $ prouve que pour tout entier $%
n\geq 1$, on a%
\begin{equation}
\binom{np-1}{p-1}\equiv 1\text{ }\func{mod}(p^{3})\text{, \ \ \ \ \ \ }p\geq
5  \label{rel74}
\end{equation}%
\begin{equation}
\binom{np-1}{p-1}\equiv 1-\frac{1}{3}n(n-1)p^{3}B_{p-3}\text{ }\func{mod}%
(p^{4})\text{, \ \ \ \ \ \ }p\geq 5.\text{\ }  \label{rel3}
\end{equation}%
Les nombres de Bernoulli $B_{n}$ \'{e}tant d\'{e}finis par leur s\'{e}rie g%
\'{e}n\'{e}ratrice 
\begin{equation*}
\frac{z}{e^{z}-1}=\sum_{n=0}^{\infty }B_{n}\frac{z^{n}}{n!}.
\end{equation*}%
En 1900 Glaisher prouve que

\begin{equation*}
\binom{2p-1}{p-1}\equiv 1+2p\sum_{k=1}^{p-1}\frac{1}{k}\text{ }\func{mod}%
(p^{4})\text{ \ \ }p\geq 3\text{.}
\end{equation*}%
En 1953, Carlitz [\cite{car1} et \cite{car2}] am\'{e}liore la congruence de
Morley en prouvant que 
\begin{equation*}
(-1)^{\frac{p-1}{2}}\binom{p-1}{\frac{p-1}{2}}\equiv 4^{p-1}+\frac{p^{3}}{12}%
B_{p-3}\func{mod}(p^{4})\text{ \ \ \ \ \ \ }p\geq 5.
\end{equation*}%
En 1995, R.J. McIntosh $\left[ \text{\cite{mac}, p. 385}\right] $ prouve que 
\begin{equation}
\binom{2p-1}{p-1}\equiv 1-p^{2}\sum_{k=1}^{p-1}\frac{1}{k^{2}}\text{ }\func{%
mod}(p^{5})\text{ \ \ \ \ \ }p\geq 7.  \label{rel4}
\end{equation}%
En 2007, Zhao \cite{zha} prouve que

\begin{equation}
\binom{2p-1}{p-1}\equiv 1+2p\sum_{k=1}^{p-1}\frac{1}{k}\text{ }\func{mod}%
(p^{5})\text{ \ \ }p\geq 7).  \label{rel4b}
\end{equation}%
En 2010, Tauraso \cite{tau} prouve que 
\begin{equation}
\binom{2p-1}{p-1}\equiv 1+2p\sum_{k=1}^{p-1}\frac{1}{k}+\frac{2}{3}%
p^{3}\sum_{k=1}^{p-1}\frac{1}{k^{3}}\text{ }\func{mod}(p^{6})\text{\ \ \ \ \
\ \ }p\geq 7.  \label{92}
\end{equation}%
et%
\begin{equation}
\binom{2p-1}{p-1}\equiv 1-2p\sum_{k=1}^{p-1}\frac{1}{k}-2p^{2}%
\sum_{k=1}^{p-1}\frac{1}{k^{2}}\text{ }\func{mod}(p^{6})\text{\ \ \ \ \ \ \ }%
p\geq 7.  \label{93}
\end{equation}%
En 2014, Me\v{s}trovi\'{c} \cite{mes} prouve que

\begin{equation}
\binom{2p-1}{p-1}\equiv 1-2p\sum_{k=1}^{p-1}\frac{1}{k}+4p^{2}\sum_{1\leq
i<j\leq p-1}\frac{1}{ij}\text{ }\func{mod}(p^{7})\text{ \ \ \ \ \ \ }p\geq 11%
\text{.}  \label{80}
\end{equation}%
La congruence \ref{80} est encore g\'{e}n\'{e}ralis\'{e}e par J. Rosen \cite%
{ros}.

\section{Enonc\'{e} du r\'{e}sultat principal}

Le th\'{e}or\`{e}me suivant constitue \`{a} la fois une g\'{e}n\'{e}%
ralisation de la congruence de Wolstenholme et de la congruence de Morley.
En exploitant la relation suivante qui d\'{e}coule du (lemme \ref{lem2})%
\begin{equation}
(-1)^{\frac{p-1}{2}}\binom{p-1}{\frac{p-1}{2}}=4^{p-1}\binom{\frac{1}{2}p-1}{%
p-1}.  \label{rel29}
\end{equation}%
Il permet aussi de retrouver toutes les nombreuses g\'{e}n\'{e}ralisations
de ces deux congruences que l'on a expos\'{e} au premier paragraphe et aussi
d'en d\'{e}couvrir d'autres.

\begin{theorem}
\label{theo}Pour tout nombre premier impair $p$ et pour tout $p$-entier $%
\alpha $, on a%
\begin{equation}
\binom{\alpha p-1}{p-1}\equiv 1-\alpha (\alpha -1)(\alpha ^{2}-\alpha
-1)p\sum_{k=1}^{p-1}\frac{1}{k}+\alpha ^{2}(\alpha -1)^{2}p^{2}\sum_{1\leq
i<j\leq p-1}\frac{1}{ij}\text{ }\func{mod}(p^{m})\text{.}  \label{rel11}
\end{equation}%
o\`{u} $m=7$ si $p\neq 7$ et $m=6$ si $p=7$.
\end{theorem}

Pour $\alpha =2$ et pour $\alpha =\frac{1}{2}$, le th\'{e}or\`{e}me \ref%
{theo} permet d'obtenir le corollaire suivant:

\begin{corollary}
\label{coro1}Pour tout nombre premier impair $p$, on a%
\begin{equation}
\binom{2p-1}{p-1}\equiv 1-2p\sum_{k=1}^{p-1}\frac{1}{k}+4p^{2}\sum_{1\leq
i<j\leq p-1}\frac{1}{ij}\text{ }\func{mod}(p^{m})  \label{rel30}
\end{equation}%
\begin{equation}
(-1)^{\frac{p-1}{2}}\binom{p-1}{\frac{p-1}{2}}\equiv 4^{p-1}\left( 1-\frac{5%
}{16}p\sum_{k=1}^{p-1}\frac{1}{k}+\frac{1}{16}p^{2}\sum_{1\leq i<j\leq p-1}%
\frac{1}{ij}\right) \text{ }\func{mod}(p^{m})  \label{rel31}
\end{equation}%
o\`{u} $m=7$ si $p\neq 7$ et $m=6$ si $p=7$.
\end{corollary}

On constate ainsi que le th\'{e}or\`{e}me \ref{theo} est bien une g\'{e}n%
\'{e}ralisation des congruences de Wolstenholme (\ref{rel70}) et de Morley(%
\ref{rel71}). En effet ces deux congruences se d\'{e}duisent respectivement
de (\ref{rel30}) et (\ref{rel31}) en observant qu'on a d'apr\`{e}s (\ref%
{rel19}) pour $m=1$ et (\ref{rel18}) pour $m=2$%
\begin{equation}
\sum_{k=1}^{p-1}\frac{1}{k}\equiv 0\func{mod}(p^{2})\text{ \ et \ }%
\sum_{1\leq i<j\leq p-1}\frac{1}{ij}\equiv 0\func{mod}(p),\text{ \ \ }p\geq 5%
\text{.}  \label{rel73}
\end{equation}%
On en d\'{e}duit aussi du th\'{e}or\`{e}me \ref{theo} et de (\ref{rel73}) le
corollaire suivant qui g\'{e}n\'{e}ralise aux $p$-entiers, la congruence de
Glaisher (\ref{rel74}).

\begin{corollary}
Pour tout $p$-entier $\alpha $, on a%
\begin{equation}
\binom{\alpha p-1}{p-1}\equiv 1\text{ }\func{mod}(p^{3}),\text{ \ \ }p\geq 5%
\text{.}  \label{rel26}
\end{equation}
\end{corollary}

On d\'{e}duit du th\'{e}or\`{e}me \ref{theo} les deux relations suivantes

\begin{equation}
\binom{2p-1}{p-1}\equiv 1-2p\sum_{k=1}^{p-1}\frac{1}{k}+4p^{2}\sum_{1\leq
i<j\leq p-1}\frac{1}{ij}\text{ }\func{mod}(p^{6})\text{,}  \label{rel32}
\end{equation}%
\begin{equation}
(-1)^{\frac{p-1}{2}}\binom{p-1}{\frac{p-1}{2}}\equiv 4^{p-1}\left( 1-\frac{5%
}{16}p\sum_{k=1}^{p-1}\frac{1}{k}+\frac{1}{16}p^{2}\sum_{1\leq i<j\leq p-1}%
\frac{1}{ij}\right) \text{ }\func{mod}(p^{6}).  \label{rel33}
\end{equation}%
Comme on a

\begin{equation*}
\sum_{1\leq i<j\leq p-1}\frac{1}{ij}=\frac{1}{2}\left( \sum_{k=1}^{p-1}\frac{%
1}{k}\right) ^{2}-\frac{1}{2}\sum_{k=1}^{p-1}\frac{1}{k^{2}},
\end{equation*}%
on en d\'{e}duit de (\ref{rel73}) que%
\begin{equation*}
p^{2}\sum_{1\leq i<j\leq p-1}\frac{1}{ij}\equiv -\frac{1}{2}%
p^{2}\sum_{k=1}^{p-1}\frac{1}{k^{2}}\func{mod}(p^{6}).
\end{equation*}%
Compte tenu de cette derni\`{e}re relation, on d\'{e}duit du th\'{e}or\`{e}%
me \ref{theo} le corollaire suivant

\begin{corollary}
Pour tout nombre premier impair $p$ et pour tout $p$-entier $\alpha $, on a%
\begin{equation}
\binom{\alpha p-1}{p-1}\equiv 1-\alpha (\alpha -1)(\alpha ^{2}-\alpha
-1)p\sum_{k=1}^{p-1}\frac{1}{k}-\frac{1}{2}\alpha ^{2}(\alpha
-1)^{2}p^{2}\sum_{k=1}^{p-1}\frac{1}{k^{2}}\text{ }\func{mod}(p^{6})
\label{rel38}
\end{equation}
\end{corollary}

Pour $\alpha =2$ et $\alpha =\frac{1}{2}$ ce corollaire nous fournit les
deux relations suivantes 
\begin{equation}
\binom{2p-1}{p-1}\equiv 1-2p\sum_{k=1}^{p-1}\frac{1}{k}-2p^{2}%
\sum_{k=1}^{p-1}\frac{1}{k^{2}}\text{ }\func{mod}(p^{6}).  \label{rel36}
\end{equation}%
\begin{equation}
(-1)^{\frac{p-1}{2}}\binom{p-1}{\frac{p-1}{2}}\equiv 4^{p-1}\left( 1-\frac{5%
}{16}p\sum_{k=1}^{p-1}\frac{1}{k}-\frac{1}{32}p^{2}\sum_{k=1}^{p-1}\frac{1}{%
k^{2}}\right) \text{ }\func{mod}(p^{6}).  \label{rel37}
\end{equation}

\bigskip De la relation (\ref{rel38}) et du lemme \ref{lem4}, on d\'{e}duit
le corollaire suivant

\begin{corollary}
\label{coro5}%
\begin{equation}
\binom{\alpha p-1}{p-1}\equiv 1-\frac{1}{3}\alpha (\alpha -1)p^{3}B_{p-3}%
\text{ }\func{mod}(p^{4})  \label{rel2}
\end{equation}
\end{corollary}

(\ref{rel2}) est une g\'{e}n\'{e}ralisation de la congruence de Glaisher (%
\ref{rel3}).

Pour $m=1$, la relation (\ref{64}) implique 
\begin{equation}
2p\sum_{k=1}^{p-1}\frac{1}{k}+p^{2}\sum_{k=1}^{p-1}\frac{1}{k^{2}}\equiv 0%
\func{mod}(p^{5})\text{ \ \ \ }p\geq 7.  \label{rel34}
\end{equation}%
On d\'{e}duit de (\ref{rel37}) et (\ref{rel38}) le corollaire suivant

\begin{corollary}
Pour tout $p$-entier $\alpha $, on a 
\begin{equation}
\binom{\alpha p-1}{p-1}\equiv 1+\alpha (\alpha -1)p\sum_{k=1}^{p-1}\frac{1}{k%
}\text{ }\func{mod}(p^{5})\text{ \ \ \ \ \ \ }p\geq 7,  \label{rel5b}
\end{equation}%
\begin{equation}
\binom{\alpha p-1}{p-1}\equiv 1-\frac{1}{2}\alpha (\alpha
-1)p^{2}\sum_{k=1}^{p-1}\frac{1}{k^{2}}\text{ }\func{mod}(p^{5})\text{ \ \ \
\ \ \ }p\geq 7\text{.}  \label{rel5}
\end{equation}
\end{corollary}

Pour $\alpha =2$, la relation (\ref{rel5}) permet de retrouver \ la
congruence de R.J. McIntosh (\ref{rel4}) et la relation (\ref{rel5b}) permet
de retrouver \ la congruence de Zhao (\ref{rel4b}) . De plus pour $\alpha =%
\frac{1}{2}$, les relations (\ref{rel36}) et (\ref{rel37}) permettent
d'obtenir le corollaire suivant

\begin{corollary}
\label{coro4} On a%
\begin{equation}
(-1)^{\frac{p-1}{2}}\binom{p-1}{\frac{p-1}{2}}\equiv 4^{p-1}\left( 1-\frac{1%
}{4}p\sum_{k=1}^{p-1}\frac{1}{k}\right) \text{ }\func{mod}(p^{5})\text{ \ \
\ \ \ \ }p\geq 7,  \label{rel6b}
\end{equation}
\begin{equation}
(-1)^{\frac{p-1}{2}}\binom{p-1}{\frac{p-1}{2}}\equiv 4^{p-1}\left( 1+\frac{1%
}{8}p^{2}\sum_{k=1}^{p-1}\frac{1}{k^{2}}\right) \text{ }\func{mod}(p^{5})%
\text{ \ \ \ \ \ }p\geq 7.  \label{rel6}
\end{equation}
\end{corollary}

A l'aide de la relation (\ref{62}) du lemme \ref{lem3}, \'{e}crite pour $m=1$%
, on d\'{e}duit que l'on a 
\begin{equation}
\sum_{k=1}^{p-1}\frac{1}{k}+\frac{1}{2}p\sum_{k=1}^{p-1}\frac{1}{k^{2}}+%
\frac{1}{6}p^{2}\sum_{k=1}^{p-1}\frac{1}{k^{3}}\equiv 0\func{mod}(p^{6})\ \
\ \ \ \ \ p\geq 11.  \label{63}
\end{equation}%
Avec (\ref{rel38}) et (\ref{63}), on d\'{e}duit le corollaire suivant.

\begin{corollary}
Pour tout nombre premier impair $p\geq 11$ et pour tout $p$-entier $\alpha $%
, on a%
\begin{equation*}
\binom{\alpha p-1}{p-1}\equiv 1+\alpha (\alpha -1)p\sum_{k=1}^{p-1}\frac{1}{k%
}+\frac{1}{6}\alpha ^{2}(\alpha -1)^{2}p^{3}\sum_{k=1}^{p-1}\frac{1}{k^{3}}%
\text{ }\func{mod}(p^{6}).
\end{equation*}
\end{corollary}

Pour $\alpha =2$ et $\alpha =\frac{1}{2}$, ce corollaire nous fournit les
deux congruences suivantes v\'{e}rifi\'{e}es pour $p\geq 11,$ 
\begin{equation*}
\binom{2p-1}{p-1}\equiv 1+2p\sum_{k=1}^{p-1}\frac{1}{k}+\frac{2}{3}%
p^{3}\sum_{k=1}^{p-1}\frac{1}{k^{3}}\text{ }\func{mod}(p^{6}),
\end{equation*}%
\begin{equation*}
(-1)^{\frac{p-1}{2}}\binom{p-1}{\frac{p-1}{2}}\equiv 4^{p-1}\left( 1-\frac{1%
}{4}p\sum_{k=1}^{p-1}\frac{1}{k}+\frac{1}{96}p^{3}\sum_{k=1}^{p-1}\frac{1}{%
k^{3}}\right) \text{ }\func{mod}(p^{6}).
\end{equation*}

Le th\'{e}or\`{e}me \ref{theo} g\'{e}n\'{e}ralise aux $p$-entiers, la
relation (\ref{80}) de Me\v{s}trovi\'{c}.

\section{Lemmes}

\begin{lemme}
\label{lem2}Pour tout entier $n\geq 1$, on a%
\begin{equation}
(-1)^{n}\binom{2n}{n}=4^{2n}\binom{n-\frac{1}{2}}{2n}  \label{rel1}
\end{equation}
\end{lemme}

\bigskip

\begin{preuve}
On a%
\begin{eqnarray*}
4^{2n}\binom{n-\frac{1}{2}}{2n} &=&\frac{4^{2n}}{(2n)!}\prod_{k=1}^{2n}(n+%
\frac{1}{2}-k) \\
&=&\frac{2^{2n}}{(2n)!}\prod_{k=1}^{2n}(2n+1-2k) \\
&=&(-1)^{n}\frac{2^{2n}}{(2n)!}\prod_{k=1}^{n}\left( 2\left( n+1-k\right)
-1\right) \prod_{k=n+1}^{2n}(2(k-n)-1),
\end{eqnarray*}%
ce qui peut s'\'{e}crire 
\begin{equation}
4^{2n}\binom{n-\frac{1}{2}}{2n}=(-1)^{n}\frac{2^{2n}}{(2n)!}\left(
\prod_{k=1}^{n}\left( 2j-1\right) \right) ^{2}.  \label{rel27}
\end{equation}%
On constate alors que 
\begin{equation}
\prod_{k=1}^{n}\left( 2j-1\right) =\frac{\prod_{j=1}^{n}(2j).%
\prod_{j=1}^{n}(2j-1)}{\prod_{j=1}^{n}(2j)}=\frac{(2n)!}{2^{n}n!}
\label{rel28}
\end{equation}%
Il r\'{e}sulte de (\ref{rel27}) et (\ref{rel28}) que 
\begin{eqnarray*}
4^{2n}\binom{n-\frac{1}{2}}{2n} &=&(-1)^{n}\frac{2^{2n}}{(2n)!}\left( \frac{%
(2n)!}{2^{n}n!}\right) ^{2} \\
&=&(-1)^{n}\frac{(2n)!}{n!n!}=(-1)^{n}\binom{2n}{n}.
\end{eqnarray*}%
En chosissant $n=\frac{p-1}{2}$ dans (\ref{rel1}), on obtient la relation (%
\ref{rel29}).
\end{preuve}

Pour tout nombre premier $p$ et pour tout entier $k$, nous d\'{e}finissons
les nombres harmoniques g\'{e}n\'{e}ralis\'{e}s $H_{m}$ par 
\begin{equation*}
H_{m}=\sum_{1\leq k_{1}<\ldots <k_{m}\leq p-1}\frac{1}{k_{1}\ldots k_{m}}%
\text{, \ \ \ \ \ \ \ \ \ \ \ \ \ \ \ \ pour }1\leq m\leq p-1
\end{equation*}%
et par convention 
\begin{equation}
H_{0}=1\text{ et }H_{m}=0\text{ pour }m\geq p.  \label{rel15}
\end{equation}%
Soit $P(x)$ le polyn\^{o}me d\'{e}fini par 
\begin{equation}
P(x)=\frac{(x-1)(x-2)\ldots (x-p+1)}{(p-1)!},  \label{rel12}
\end{equation}%
On a alors 
\begin{eqnarray*}
P(x) &=&(-1)^{p-1}\prod_{k=1}^{p-1}\frac{k-x}{k} \\
&=&\prod_{i=1}^{p-1}\left( 1-\frac{x}{k}\right) .
\end{eqnarray*}%
On en d\'{e}duit que 
\begin{equation}
P(x)=\sum_{k=0}^{p-1}\left( -1\right) ^{k}H_{k}x^{k}.  \label{rel13}
\end{equation}

La preuve du th\'{e}or\`{e}me principal repose essentiellement sur le lemme
suivant

\begin{lemme}
Pour tout nombre premier $p$ impair et pour tout entier $m\geq 1$, on a

\begin{enumerate}
\item 
\begin{equation}
H_{2m-1}-mpH_{2m}=\frac{1}{2}p^{2}\sum_{k=2m+1}^{p-1}(-1)^{k}\binom{k}{2m-1}%
\text{ }p^{k-2m-1}H_{k},  \label{rel17}
\end{equation}

\item 
\begin{equation}
H_{m}\equiv 0\text{ }\func{mod}(p)\text{, pour }m\neq p-1,  \label{rel18}
\end{equation}

\item 
\begin{equation}
H_{m}\equiv 0\text{ }\func{mod}(p^{2})\text{, pour }m\text{ impair et }m\neq
p-2,  \label{rel19}
\end{equation}

\item 
\begin{equation}
H_{2m-1}-mpH_{2m}\equiv 0\text{ }\func{mod}(p^{4})\text{, pour }2m+1\neq p-2.
\label{rel20}
\end{equation}
\end{enumerate}
\end{lemme}

\begin{preuve}
\begin{enumerate}
\item La relation (\ref{rel12}) nous permet de constater que l'on a 
\begin{equation*}
P(x)=P(p-x),
\end{equation*}%
ce qui peut s'\'{e}crire, en exploitant la relation (\ref{rel13}) et la
convention (\ref{rel15}):%
\begin{equation}
\sum_{k\geq 0}\left( -1\right) ^{k}H_{k}x^{k}=\sum_{k\geq 0}\left( -1\right)
^{k}H_{k}(p-x)^{k}.  \label{rel14}
\end{equation}%
En identifiant coefficients de $x^{2m-1}$ dans chacun des deux membres de (%
\ref{rel14}), on obtient%
\begin{eqnarray*}
-H_{2m-1} &=&\sum_{k\geq 0}\left( -1\right) ^{k}H_{k}\binom{k}{2m-1}%
p^{k-2m+1}(-1)^{2m-1} \\
&=&-\sum_{k=2m-1}^{p-1}\left( -1\right) ^{k}\binom{k}{2m-1}p^{k-2m+1}H_{k} \\
&=&H_{2m-1}-2mpH_{2m}-p^{2}\sum_{k=2m+1}^{p-1}\left( -1\right) ^{k}\binom{k}{%
2m-1}p^{k-2m-1}H_{k}.
\end{eqnarray*}%
La relation (\ref{rel17}) en r\'{e}sulte. Remarquons qu'on d\'{e}duit de
cette relation que $H_{2m-1}\equiv 0$ $\func{mod}(p)$ pour tout $m\geq 1$.

\item Il s'agit d'un r\'{e}sulat bien connu qu'on peut d\'{e}duire \
facilement du fait que le polyn\^{o}me $P(x)$ consid\'{e}r\'{e} comme un
polyn\^{o}me \`{a} coefficients dans le corps $%
\mathbb{Z}
/p%
\mathbb{Z}
$ s'\'{e}crit gr\^{a}ce au petit th\'{e}or\`{e}me de Fermat $P(x)=1-x^{p-1}$%
. On en d\'{e}duit aussi que $H_{p-1}=\frac{1}{(p-1)!}\equiv -1$ $\func{mod}%
(p)$.

\item D'apr\`{e}s la relation (\ref{rel17}), on a 
\begin{equation}
H_{2m-1}\equiv mpH_{2m}\text{ }\func{mod}(p^{2}).  \label{rel21}
\end{equation}%
On a alors $H_{2m-1}\equiv 0$ $\func{mod}(p^{2})$ pour $2m-1\neq p-2$ car on
a $H_{2m}\equiv 0$ $\func{mod}(p)$ d'apr\`{e}s (\ref{rel18}). Remarquons que
si $2m-1=p-2$, on a $H_{2m-1}=H_{p-2}\equiv \frac{p-1}{2}pH_{p-1}$ $\equiv 
\frac{p}{2}\func{mod}(p^{2})$.

\item D'apr\`{e}s la relation (\ref{rel17}), on a%
\begin{equation}
H_{2m-1}-mpH_{2m}\equiv -\frac{1}{2}p^{2}\binom{2m+1}{2}H_{2m+1}+\frac{1}{2}%
p^{3}\binom{2m+2}{3}H_{2m+2}\text{ }\func{mod}(p^{4}).  \label{rel22}
\end{equation}%
Si $2m+1\neq p-2$, on a alors $2m+2\neq p-1$ et on d\'{e}duit de (\ref{rel19}%
) et (\ref{rel18}) que $H_{2m+1}\equiv 0$ $\func{mod}(p^{2})$ et $%
H_{2m+2}\equiv 0$ $\func{mod}(p)$. En tenant compte de ces deux derni\`{e}%
res congruences dans (\ref{rel22}), la relation (\ref{rel20}) en r\'{e}sulte.
\end{enumerate}
\end{preuve}

\begin{remarque}
Il est facile de prouver \`{a} l'aide de la relation (\ref{rel22}) que si $%
2m+1=p-2$, alors $.$ 
\begin{equation*}
H_{2m-1}-mpH_{2m}=H_{p-4}-\left( \frac{p-3}{2}\right) pH_{p-3}\equiv -\frac{%
p^{3}}{4}\text{ }\func{mod}(p^{4})\text{. }
\end{equation*}%
Ainsi, pour tout $m\geq 1$, on a $H_{2m-1}-mpH_{2m}\equiv 0$ $\func{mod}%
(p^{3})$.
\end{remarque}

\begin{lemme}
\label{lem3} Pour tout entier $m\geq 1$, on a

\begin{enumerate}
\item 
\begin{equation*}
\sum_{k=1}^{p-1}\frac{1}{k^{m}}\equiv \left\{ 
\begin{array}{ll}
\text{ \ }0\text{ }\func{mod}(p)\ \ \ \ \  & \text{si }p-1\nmid m \\ 
-1\text{ }\func{mod}(p) & \text{si }p-1\mid m%
\end{array}%
\right. .
\end{equation*}

\item 
\begin{equation*}
\sum_{k=1}^{p-1}\frac{1}{k^{m}}\equiv \left\{ 
\begin{array}{ll}
\text{ \ \ }0\text{ \ \ }\func{mod}(p^{2})\ \ \ \ \  & \text{si }m\text{ est
impair et }p-1\nmid m+1 \\ 
\frac{1}{2}mp\text{ }\func{mod}(p^{2}) & \text{si }m\text{ est impair et }%
p-1\mid m+1%
\end{array}%
\right. .
\end{equation*}

\item Pour $m$ impair, on a%
\begin{equation}
2\sum_{k=1}^{p-1}\frac{1}{k^{m}}+mp\sum_{k=1}^{p-1}\frac{1}{k^{m+1}}\equiv
\left\{ 
\begin{array}{ll}
\equiv 0\text{ \ \ \ \ \ \ \ \ \ \ \ \ \ \ \ \ \ \ \ \ \ \ \ \ \ \ \ \ \ \ \ 
}\func{mod}(p^{3})\ \  &  \\ 
\equiv 0\text{ \ \ \ \ \ \ \ \ \ \ \ \ \ \ \ \ \ \ \ \ \ \ \ \ \ \ \ \ \ \ \ 
}\func{mod}(p^{4}) & \text{si}\ p-1\nmid m+3 \\ 
\equiv -\frac{1}{12}m(m+1)(m+2)p^{3}\text{ }\func{mod}(p^{4}) & \text{si}\
p-1\mid m+3%
\end{array}%
\right.  \label{64}
\end{equation}

\item Pour $m$ impair, on a%
\begin{equation}
\sum_{k=1}^{p-1}\frac{1}{k^{m}}+\frac{1}{2}mp\sum_{k=1}^{p-1}\frac{1}{k^{m+1}%
}+\frac{m(m+1)}{12}p^{2}\sum_{k=1}^{p-1}\frac{1}{k^{m+2}}\equiv 0\func{mod}%
(p^{6})\ \ \ \ \ \text{si}\ p-1\nmid m+5\text{.}\ \   \label{62}
\end{equation}
\end{enumerate}
\end{lemme}

\begin{preuve}
\begin{enumerate}
\item Soit $g\in 
\mathbb{Z}
$ tel que $\overline{g}$ soit un g\'{e}n\'{e}rateur du groupe cyclique $%
\left( 
\mathbb{Z}
/p%
\mathbb{Z}
\right) ^{\ast }$. L'application $x\rightarrow x^{-1}$ de $\left( 
\mathbb{Z}
/p%
\mathbb{Z}
\right) ^{\ast }$ dans lui m\^{e}me \'{e}tant bijective, on a

\begin{equation*}
\sum_{k=1}^{p-1}\frac{1}{k^{m}}\equiv \sum_{k=1}^{p-1}k^{m}\equiv
\sum_{j=0}^{p-2}\left( g^{j}\right) ^{m}\equiv \sum_{j=0}^{p-2}\left(
g^{m}\right) ^{j}\equiv 0\func{mod}(p).
\end{equation*}

On a alors 
\begin{equation*}
\left( g^{m}-1\right) \sum_{k=1}^{p-1}\frac{1}{k^{m}}\equiv \left(
g^{m}-1\right) \sum_{j=0}^{p-2}\left( g^{m}\right) ^{j}=\left( g^{m}\right)
^{p-1}-1=\left( g^{p-1}\right) ^{m}-1\equiv 0\func{mod}(p)\text{.}
\end{equation*}

On en d\'{e}duit que si $p-1\nmid m$, on a $g^{m}-1\not\equiv 0\func{mod}(p)$
et $\sum_{k=1}^{p-1}\frac{1}{k^{m}}\equiv 0\func{mod}(p)$.

Si $p-1\mid m$, on a $\frac{1}{k^{m}}\equiv 1\func{mod}(p)$ pour $1\leq
k\leq p-1$ et $\sum_{k=1}^{p-1}\frac{1}{k^{m}}\equiv
\sum_{k=1}^{p-1}1=p-1\equiv -1\func{mod}(p)$.

\item Pour $m$ impair, on a%
\begin{eqnarray}
\sum_{k=1}^{p-1}\frac{1}{k^{m}} &=&\frac{1}{2}\sum_{k=1}^{p-1}\frac{1}{k^{m}}%
+\frac{1}{2}\sum_{k=1}^{p-1}\frac{1}{(p-k)^{m}}  \notag \\
&=&\frac{1}{2}\sum_{k=1}^{p-1}\frac{(p-k)^{m}+k^{m}}{k^{m}(p-k)^{m}}\equiv -%
\frac{1}{2}mp\sum_{k=1}^{p-1}\frac{1}{k^{m+1}}\func{mod}(p^{2}).  \label{61}
\end{eqnarray}%
Si $p-1\nmid m+1$, on a $\sum_{k=1}^{p-1}\frac{1}{k^{m+1}}\equiv 0\func{mod}%
(p)$ et (\ref{61}) implique $\sum_{k=1}^{p-1}\frac{1}{k^{m}}\equiv 0\func{mod%
}(p^{2})$.

Si $\ p-1\mid m+1$; alors $\sum_{k=1}^{p-1}\frac{1}{k^{m+1}}\equiv -1\func{%
mod}(p)$ et (\ref{61}) implique $\sum_{k=1}^{p-1}\frac{1}{k^{m}}\equiv \frac{%
1}{2}mp$ $\func{mod}(p^{2})$.

\item On a pour $1\leq k\leq p-1$%
\begin{equation*}
\frac{1}{\left( 1-\frac{p}{k}\right) ^{m}}\equiv 1+m\frac{p}{k}+\frac{m(m+1)%
}{2}\frac{p^{2}}{k^{2}}+\frac{m(m+1)(m+2)}{6}\frac{p^{3}}{k^{3}}\text{ }%
\func{mod}(p^{4})\text{.}
\end{equation*}%
On en d\'{e}duit que pour $m$ impair, on a 
\begin{eqnarray*}
\sum_{k=1}^{p-1}\frac{1}{(p-k)^{m}} &=&-\sum_{k=1}^{p-1}\frac{1}{k^{m}(1-%
\frac{p}{k})^{m}} \\
&\equiv &\sum_{k=1}^{p-1}\left( -\frac{1}{k^{m}}-m\frac{p}{k^{m+1}}-\frac{%
m(m+1)}{2}\frac{p^{2}}{k^{m+2}}-\frac{m(m+1)(m+2)}{6}\frac{p^{3}}{k^{m+3}}%
\right) \text{ }\func{mod}(p^{4}).
\end{eqnarray*}%
Il en r\'{e}sulte que%
\begin{eqnarray}
2\sum_{k=1}^{p-1}\frac{1}{k^{m}} &=&\sum_{k=1}^{p-1}\frac{1}{k^{m}}%
+\sum_{k=1}^{p-1}\frac{1}{(p-k)^{m}}  \notag \\
&\equiv &-mp\sum_{k=1}^{p-1}\frac{1}{k^{m+1}}-\frac{m(m+1)}{2}%
p^{2}\sum_{k=1}^{p-1}\frac{1}{k^{m+2}}-\frac{m(m+1)(m+2)}{6}%
p^{3}\sum_{k=1}^{p-1}\frac{1}{k^{m+3}}\text{ }\func{mod}(p^{4}).  \label{60}
\end{eqnarray}%
Or $m+2$ est \ impair. On a donc $p-1\nmid m+2$ et $\sum_{k=1}^{p-1}\frac{1}{%
k^{m+2}}\equiv 0\func{mod}(p)$. On d\'{e}duit alors de (\ref{60}) que 
\begin{equation*}
2\sum_{k=1}^{p-1}\frac{1}{k^{m}}+mp\sum_{k=1}^{p-1}\frac{1}{k^{m+1}}\equiv 0%
\text{ }\func{mod}(p^{3}).
\end{equation*}%
Si $p-1\nmid m+3$, on a \`{a} la fois $\sum_{k=1}^{p-1}\frac{1}{k^{m+2}}%
\equiv 0$ $\func{mod}(p^{2})$ et $\sum_{k=1}^{p-1}\frac{1}{k^{m+3}}\equiv 0%
\func{mod}(p)$. Dans ce cas, on d\'{e}duit de (\ref{60}) que 
\begin{equation*}
2\sum_{k=1}^{p-1}\frac{1}{k^{m}}\equiv -mp\sum_{k=1}^{p-1}\frac{1}{k^{m+1}}%
\text{ }\func{mod}(p^{4}).
\end{equation*}%
Si$\ p-1\mid m+3$, on a $\sum_{k=1}^{p-1}\frac{1}{k^{m+2}}\equiv \frac{1}{2}%
(m+2)p$ $\ \func{mod}(p^{2})$ et $\sum_{k=1}^{p-1}\frac{1}{k^{m+3}}\equiv -1$
$\ \func{mod}(p)$. Dans ce cas, on d\'{e}duit de (\ref{60}) que%
\begin{eqnarray*}
2\sum_{k=1}^{p-1}\frac{1}{k^{m}}+mp\sum_{k=1}^{p-1}\frac{1}{k^{m+1}} &\equiv
&-\frac{m(m+1)(m+2)}{4}p^{3}+\frac{m(m+1)(m+2)}{6}p^{3}\text{ } \\
&\equiv &-\frac{1}{12}m(m+1)(m+2)\text{ }\func{mod}(p^{4}).
\end{eqnarray*}

\item Si$\ $\ $m$ impair et si $p-1\nmid m+5$, on a 
\begin{equation}
2\sum_{k=1}^{p-1}\frac{1}{k^{m}}\equiv -mp\sum_{k=1}^{p-1}\frac{1}{k^{m+1}}-%
\frac{m(m+1)}{2}p^{2}\sum_{k=1}^{p-1}\frac{1}{k^{m+2}}-\frac{m(m+1)(m+2)}{6}%
p^{3}\sum_{k=1}^{p-1}\frac{1}{k^{m+3}}\text{ }\func{mod}(p^{4}).  \label{65}
\end{equation}%
On a aussi d'apr\`{e}s (\ref{64}) 
\begin{equation}
2p^{2}\sum_{k=1}^{p-1}\frac{1}{k^{m+2}}+(m+2)p^{3}\sum_{k=1}^{p-1}\frac{1}{%
k^{m+3}}\equiv 0\text{ }\func{mod}(p^{6}).  \label{66}
\end{equation}%
On d\'{e}duit de (\ref{65}) et (\ref{66})%
\begin{equation*}
2\sum_{k=1}^{p-1}\frac{1}{k^{m}}\equiv -mp\sum_{k=1}^{p-1}\frac{1}{k^{m+1}}-%
\frac{m(m+1)}{2}p^{2}\sum_{k=1}^{p-1}\frac{1}{k^{m+2}}+\frac{m(m+1)}{3}\text{
}p^{2}\sum_{k=1}^{p-1}\frac{1}{k^{m+2}}\func{mod}(p^{4}).
\end{equation*}%
La relation (\ref{62}) en r\'{e}sulte.
\end{enumerate}
\end{preuve}

\begin{lemme}
\label{lem4}On a%
\begin{eqnarray}
\sum_{k=1}^{p-1}\frac{1}{k} &\equiv &-\frac{1}{3}p^{2}B_{p-3}\func{mod}%
(p^{3}),\text{ \ \ \ \ \ }p\geq 5  \label{90} \\
\sum_{k=1}^{p-1}\frac{1}{k^{2}} &\equiv &\frac{2}{3}pB_{p-3}\func{mod}%
(p^{2}),\text{ \ \ \ \ \ \ \ \ \ }p\geq 5.  \label{91}
\end{eqnarray}
\end{lemme}

On trouvera dans \cite{ji} une preuve tr\'{e}s d\'{e}taill\'{e}e de la
relation (\ref{90}) qui est un r\'{e}sultat d\^{u} \`{a} Glaisher \cite{gla1}%
. La realtion (\ref{91}) se d\'{e}duit de (\ref{90}) et de (\ref{64}) \'{e}%
crite pour $m=1.$

\section{Preuve du th\'{e}or\`{e}me \protect\ref{theo}}

D'apr\`{e}s la relation (\ref{rel12}), on a 
\begin{equation*}
\binom{\alpha p-1}{p-1}=P(\alpha p)=\sum_{k=0}^{p-1}\left( -\alpha \right)
^{k}H_{k}p^{k}.
\end{equation*}%
On en d\'{e}duit que%
\begin{equation}
\binom{\alpha p-1}{p-1}=\sum_{k=0}^{4}\left( -\alpha \right)
^{k}H_{k}p^{k}+\left( -\alpha \right) ^{5}H_{5}p^{5}+\left( -\alpha \right)
^{6}H_{6}p^{6}\func{mod}(p^{7}).  \label{rel23}
\end{equation}%
Or, d'apr\`{e}s les relations (\ref{rel18}) et (\ref{rel19}), on a%
\begin{equation}
\left( -\alpha \right) ^{5}H_{5}p^{5}\equiv 0\text{ }\func{mod}(p^{7})\text{
\ et \ }\left( -\alpha \right) ^{6}H_{6}p^{6}\equiv 0\text{ }\func{mod}%
(p^{7}),  \label{rel24}
\end{equation}%
pourvu que $5\neq p-2$ et $6\neq p-1$, c'est \`{a} dire $p\neq 7$. Il suffit
donc de choisir $p\geq 11$ pour r\'{e}aliser ces conditions.

Ainsi pour $p\geq 11$, Il r\'{e}sulte des relations (\ref{rel23}) et (\ref%
{rel24}) que l'on a%
\begin{equation}
\binom{\alpha p-1}{p-1}\equiv \sum_{k=0}^{4}\left( -\alpha \right)
^{k}H_{k}p^{k}\func{mod}(p^{7}).  \label{rel7}
\end{equation}

Pour $\alpha =1$, nous d\'{e}duisons de (\ref{rel7}) la relation%
\begin{equation}
\sum_{k=1}^{4}\left( -\alpha \right) ^{k}H_{k}p^{k}\equiv 0\func{mod}(p^{7}).
\label{rel8}
\end{equation}%
D'autre part, du fait que $p\geq 11$, on a d'apr\`{e}s (\ref{rel20}) du
lemme, 
\begin{equation}
p^{3}H_{3}-2p^{4}H_{4}\equiv 0\text{ }\func{mod}(p^{7}).  \label{rel9}
\end{equation}

Des relations (\ref{rel7}), (\ref{rel8}) et (\ref{rel9}), on d\'{e}duit que
pour tous $p$-entiers $\lambda $ et $\mu $, on a 
\begin{equation}
\binom{\alpha p-1}{p-1}\equiv \sum_{k=0}^{4}\left( -\alpha \right)
^{k}H_{k}p^{k}+\lambda \left( \sum_{k=1}^{4}\left( -\alpha \right)
^{k}H_{k}p^{k}\right) +\mu \left( p^{3}H_{3}-2p^{4}H_{4}\right) \text{ }%
\func{mod}(p^{7}).  \label{rel10}
\end{equation}%
Autrement dit, on a 
\begin{equation*}
\binom{\alpha p-1}{p-1}\equiv \sum_{k=0}^{4}A_{k}H_{k}p^{k}\func{mod}(p^{7}),
\end{equation*}%
avec%
\begin{eqnarray*}
A_{0} &=&1, \\
A_{1} &=&-\alpha -\lambda , \\
A_{2} &=&\alpha ^{2}+\lambda , \\
A_{3} &=&-\alpha ^{3}-\lambda +11, \\
A_{4} &=&\alpha ^{4}+\lambda -2\mu .
\end{eqnarray*}%
Choisissons $\lambda $ et $\mu $ tels que $A_{3}=A_{4}=0$, on obtient 
\begin{equation*}
\lambda =\alpha ^{4}-2\alpha ^{3}\text{ \ \ \ et \ \ \ }\mu =\alpha
^{4}-\alpha ^{3}.
\end{equation*}%
Avec ce choix de $\lambda $ et $\mu $, la relation (\ref{rel10}) devient 
\begin{equation*}
\binom{\alpha p-1}{p-1}\equiv 1-(\alpha ^{4}-2\alpha ^{3}+\alpha
)pH_{1}+(\alpha ^{4}-2\alpha ^{3}+\alpha ^{2})p^{2}H_{2}\func{mod}(p^{7}).
\end{equation*}%
Ce qui nous fournit bien la relation (\ref{rel11}).

Si $p=7$, un calcul direct donne%
\begin{equation*}
\binom{\alpha p-1}{p-1}-\left( 1-\alpha (\alpha -1)(\alpha ^{2}-\alpha
-1)p\sum_{k=1}^{p-1}\frac{1}{k}+\alpha ^{2}(\alpha -1)^{2}p^{2}\sum_{1\leq
i<j\leq p-1}\frac{1}{ij}\right) =\frac{\alpha ^{3}(\alpha -1)^{3}}{720}%
7^{6}\equiv 0\func{mod}(7^{6})
\end{equation*}%
et permet de conclure. La preuve du th\'{e}or\`{e}me est compl\`{e}te.


\begin{thebibliography}{99}
\bibitem{bab} C. Babbage, Demonstration of a theorem relating to prime
numbers, Edinburgh Philosophical J. 1 (1819), 46--49.

\bibitem{car1} L. Carlitz, A Theorem of Glaisher, Canadian Journal of
Mathematics \textbf{5} (1953):306-316.

\bibitem{car2} L. Carlitz, Note on a Theorem of Glaisher. Journal of the
London Mathematical Society. \textbf{28} (1953): 245-246.

\bibitem{gla} J. W. L. Glaisher, Congruences relating to the sums of
products of the first n numbers and to other sums and products, Quart. J.
Math. 31 (1899), 2--35.

\bibitem{gla1} J.W.L. Glaisher, Congruences relating to the sums of products
of the first n numbers and to other sums of products, Q. J. Math. 31 (1900),
1--35.

\bibitem{gla2} J.W.L. Glaisher, On the residues of the sums of products of
the first $p-1$ numbers,and their powers, to modulus $p^{2}$ or $p^{3}$, Q.
J. Math. 31 (1900), 321--353.

\bibitem{har} G.H. Hardy, E.M. Wright, An Introduction to the Theory of
Numbers, Clarendon, Oxford, 1980.

\bibitem{ji} C. G. Ji, A simple proof of a curious congruence by Zhao, Proc.
Amer. Math. Soc., 133(2005): 3469-3472.

\bibitem{leh} E. Lehmer, On congruences involving Bernoulli numbers and the
quotients of Fermat and Wilson, Ann. of Math. (2) 39 (1938) 350--360.

\bibitem{mac} R.J. McIntosh, On the converse of Wolstenholme's Theorem, Acta
Arith. 71 (1995), 381--389.

\bibitem{mes} R. Me\v{s}trovi\'{c}, On the mod $p^{7}$ determination of $%
\binom{2p-1}{p-1}$, Rocky Mountain Journal of Mathematics, 44 (2), (2014),
633-648.

\bibitem{mes2} R. Me\v{s}trovi\'{c}, Wolstenholme's theorem: its
generalizations and extensions in the last hundred and fifty years $(1862$-$%
2012)$, eprint: arXiv:1111.3057 v2.

\bibitem{mor} F. Morley, Note on the congruence $2^{4n}$ $\equiv $ $(-)^{n}$(%
$2n)!/(n!)^{2},$ where $2n+1$ is a prime, Ann. of Math. 9 (1894/95), no.
1-6, 168--170.

\bibitem{ros} J. Rosen. Multiple harmonic sums and Wolstenholme's theorem.
International Journal of Number Theory, 9(8):2033--2052, 2013.

\bibitem{tau} R. Tauraso, More congruences for central binomial
coefficients, J. Number Theory \textbf{130} (2010), 2639--2649.

\bibitem{wol} J. Wolstenholme, On certain properties of prime numbers, Quart
J. Math. \textbf{5} (1862) 35-39.

\bibitem{zha} J. Zhao, Bernoulli Numbers, Wolstenholme's theorem, and $p^{5}$
variations of Lucas'theorem, J. Number Theory 123 (2007), 18--26.
\end{thebibliography}
\end{document}